\documentclass[a4paper,11pt]{article}
\usepackage{amsmath}
\usepackage{amssymb}
\usepackage{amsthm}
\usepackage[usenames]{color}
\usepackage{graphicx}
\usepackage[american]{babel}
\usepackage{dsfont}
\usepackage{hyperref}
\hypersetup{
  pdftitle={A variational principle for hyperideal circle patterns},
  pdfauthor={Boris A. Springborn},
  pdfsubject={},
  pdfkeywords={},
  pdfpagemode=None,
  pdfstartview=FitH,
  pdfnewwindow=true,
  breaklinks=true,
  pagebackref=true,
  colorlinks=true,
  linkcolor=black,
  anchorcolor=black,
  citecolor=black,
  filecolor=black,
  menucolor=black,
  pagecolor=black,
  urlcolor=blue,
  bookmarksopen=false,
}

\title{A discrete Laplace-Beltrami operator for simplicial surfaces}
\author{ Alexander I.~Bobenko \and Boris A.~Springborn}
\date{February 24, 2006}

\newcommand{\R}{\mathds R}
\newcommand{\Z}{\mathds Z}

\newcommand{\hrm}{\operatorname{\mathit{hrm}}}
\newcommand{\conv}{\operatorname{conv}}
\newcommand{\interior}{\operatorname{int}}
\newcommand{\power}{\operatorname{pow}}
\newcommand{\ie}{{\it i.e.}}
\newcommand{\etal}{{\it et~al.}}
\newcommand{\eg}{{\it e.g.}}
\newcommand{\apriori}{{\it a priori}}

\newtheorem{theorem}{Theorem}
\newtheorem*{rippastheorem}{Rippa's Theorem}
\newtheorem{lemma}[theorem]{Lemma}
\newtheorem{proposition}[theorem]{Proposition}

\newtheorem{definition}[theorem]{Definition}
\theoremstyle{remark}
\newtheorem*{remark}{Remark}

\begin{document}
\maketitle {\renewcommand{\thefootnote}{} \footnote[0]{Research for this
    article was supported by the DFG Research Unit 565 ``Polyhedral
    Surfaces'' and the DFG Research Center \textsc{Matheon} ``Mathematics for
    key technologies'' in Berlin.}}
\begin{center}
  Institut f\"ur Mathematik,
  Technische Universit\"at Berlin,\\
  Strasse~des 17.~Juni 136, 10623 Berlin, Germany\\
  \medskip {\tt bobenko@math.tu-berlin.de}\\
  {\tt springb@math.tu-berlin.de}
\end{center}

%--------------------------------------------------------------------------
\begin{abstract}
  We define a discrete Laplace-Beltrami operator for simplicial surfaces
  (Definition~\ref{d.Laplace}). It depends only on the intrinsic geometry of
  the surface and its edge weights are positive. Our Laplace operator is
  similar to the well known finite-elements Laplacian (the so called ``cotan
  formula'') except that it is based on the intrinsic Delaunay triangulation
  of the simplicial surface. This leads to new definitions of discrete
  harmonic functions, discrete mean curvature, and discrete minimal surfaces.
  The definition of the discrete Laplace-Beltrami operator depends on the
  existence and uniqueness of Delaunay tessellations in piecewise flat
  surfaces. While the existence is known, we prove the uniqueness. Using
  Rippa's Theorem we show that, as claimed, Musin's harmonic index provides
  an optimality criterion for Delaunay triangulations, and this can be
  used to prove that the edge flipping algorithm terminates also in the
  setting of piecewise flat surfaces.
\end{abstract}


\section{Dirichlet energy of piecewise linear
functions} \label{s.energy}

\begin{figure}[htbp]
\begin{center}
\parbox[l]{0.2\textwidth}{
    \makebox[0pt][l]{\hspace{-10pt}\input{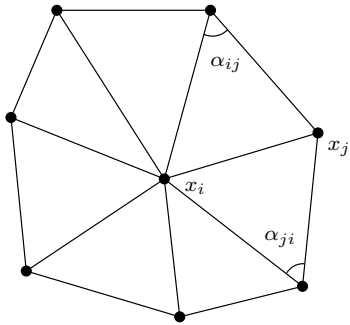}}}
    \hfill
\end{center}
\caption{The $\alpha$-angles of an internal edge.}
\label{f.angles}
\end{figure}

Let $\mathcal S$ be a \emph{simplicial surface}\/ in 3-dimensional Euclidean
space, \ie\ a geometric simplicial complex in $\R^3$ whose carrier $S$\/ is a
$2$-dimensional submanifold, possibly with boundary. We assume $\mathcal S$
to be finite.  Let $V=\{x_1,\ldots,x_{|V|}\}$, $E$, and $F$ be the sets of
vertices, edges and (triangular) faces of $\mathcal S$. Let
$f:S\rightarrow\R$ be a piecewise linear (PL) function on $S$ (linear on each
simplex of $\mathcal S$). Then the gradient $\nabla f$ is constant on each
triangle. The \emph{Dirichlet energy}\/ $ E(f) = \tfrac{1}{2}\int_S \|\nabla
f\|^2 $ is
\begin{equation*}
  E(f)=\tfrac{1}{2}\sum_{(x_i,x_j)\in E} w_{ij}\,(f(x_i)-f(x_j))^2,
\end{equation*}
where the edge weights are
\begin{equation*}
  w_{ij}=
  \left\{
    \begin{array}{ll}
      \frac{1}{2}(\cot\alpha_{ij}+\cot\alpha_{ji}) & \text{for interior edges}\\
      \frac{1}{2}\cot\alpha_{ij} &\text{for boundary edges}
    \end{array}
  \right.
\end{equation*}
and $\alpha_{ij}$, $\alpha_{ji}$ are the angle(s) opposite edge $(x_i,x_j)$
in the adjacent triangle(s) (see Figure~\ref{f.angles}). This formula was, it
seems, first derived by Duffin~\cite{duffin1}, who considers triangulated
planar regions. It follows (by summation over the triangles) from the
observation that the Dirichlet energy of a linear function on a triangle
$(x_1,x_2,x_3)$ with angles $\alpha_1$, $\alpha_2$, $\alpha_3$ is
$E(f_{\mid(x_1,x_2,x_3)})=
\frac{1}{4}\sum_{i\in\Z/3\Z}\cot\alpha_{i}\,(f(x_{i+1})-f(x_{i+2}))^2$.

In analogy to the smooth case, the Laplace operator is defined as the
gradient of the Dirichlet energy. (We identify the vector space of PL
functions $S\rightarrow\R$ with the the vector space $\R^V$ of functions on
the vertices.) By differentiating $E(f)$ with respect to the value of $f$ at
a vertex $x_i\in V$ one obtains the ``cotan formula'' for the Laplace
operator:
\begin{equation*}%                            \label{eq.Laplace_PP}
\Delta f(x_i)=\sum_{x_j\in V:(x_i,x_j)\in E}
w_{ij}\,(f(x_i)-f(x_j)).
\end{equation*}
Dziuk was the first to treat a finite element approach for the Laplace
operator on simplicial surfaces, but without stating the cotan formula
explicitely~\cite{dziuk}. It seems to have been rediscovered by Pinkall and
Polthier in their investigation of discrete minimal surfaces~\cite{PP}, and
turned out to be extremely important in geometry processing where it found
numerous applications, \eg~\cite{DMA}, \cite{BK} to name but two. In
particular, harmonic parameterizations $u:V\to {\R}^2$ are used in computer
graphics for texture mapping. The cotan-formula also forms the basis for a
theory of discrete holomorphic functions and discrete Riemann
surfaces~\cite{duffin2} \cite{Me}.

Two important \emph{disadvantages}\/ of this definition of a discrete Laplace
operator are:

\emph{1.~The weights may be negative.} The properties of the discrete Laplace
operator with {\em positive weights} ($w_{ij}>0$) are analogous to the
properties of the classical Laplace-Beltrami operator on a surface with
Riemannian metric. In particular the maximum principle holds. But some
weights $w_{ij}$ may be negative, and this leads to unpleasant phenomena: The
maximum principle does not hold. As a consequence, a vertex of a discrete
minimal surface (as defined by Pinkall \& Polthier~\cite{PP}) may not be
contained in the convex hull of its neighbors. In texture mapping
applications negative weights are undesirable because they lead to ``flipped
triangles''. In practice various tricks are used to avoid negative weights.

\emph{2.~The definition is not purely intrinsic.} The classical
Laplace-Beltrami operator is {\em intrinsic}\/ to a Riemannian manifold: It
depends only on the Riemannian metric. This is not the case with the discrete
Laplace operator defined above. Two simplicial surfaces which are isometric
but which are not triangulated in the same way give in general rise to
different Laplace operators.  As the simplest example, consider the two
triangulations of a planar quadrilateral. They lead to different discrete
Laplace operators.  (Planarity is not what causes the problem since the
quadrilateral may also be folded along either of its diagonals.)

The key idea of this paper is that one can avoid both the above shortcomings
by using the \emph{intrinsic} Delaunay triangulation of the surface $S$ to
define the discrete Laplace operator (Definition~\ref{d.Laplace}) instead of
the triangulation that comes from the simplicial complex $\mathcal S$.

\section{Delaunay triangulations of piecewise flat surfaces}
\label{s.Delaunay}

This section provides the necessary background on Delaunay tessellations of
piecewise flat surfaces. We decided to give a detailed exposition because not
all necessary proofs can be found elsewhere.

The concept of a Delaunay triangulation in $n$-dimensional Euclidean space
goes back to Delaunay~\cite{Del}. Piecewise flat surfaces
(Definition~\ref{def:PF_surface}) were studied by (his student)
Alexandrov~\cite{alexandrov} and more recently by Troyanov~\cite{Tro}. The
idea of a Delaunay triangulation of a piecewise flat surface was apparently
first considered by Rivin~\cite[Sec.~10]{Riv}. The vertex set of the Delaunay
triangulation is assumed to contain the set of cone points of the piecewise
flat surface, so that the surface is flat away from the vertices. (This is
very different from considering Delaunay triangulations in surfaces with
Riemannian metric~\cite{leibon}.) Rivin claims but does not prove an
existence and uniqueness theorem for Delaunay triangulations in piecewise
flat surfaces. His proof that the edge flipping algorithm terminates is
flawed (see the discussion after Proposition~\ref{t.flipalgorithm} below). A
correct proof was given by Indermitte \etal~\cite{ILTC}. (They seem to miss a
small detail, a topological obstruction to edge-flipability. See our proof
of Proposition~\ref{prop:flippable}.)  Furthermore, for our definition of the
discrete Laplace-Beltrami operator we also need the \emph{uniqueness} of the
Delaunay tessellation, and this question has not been addressed properly.
Rivin and Indermitte~\etal\ define Delaunay triangulations by the local
Delaunay criterion (see Definition~\ref{def:locally_delaunay}), and infer
existence via the edge flipping algorithm (see
Proposition~\ref{t.flipalgorithm}). To obtain uniqueness, we will define the
Delaunay tessellation (whose faces are generically but not always triangular)
via a global empty circle criterion
(Definition~\ref{def:delaunay_tesselation}). In a piecewise flat surface, the
``empty circumcircles'' are immersed empty disks
(Definition~\ref{def:immersed_empty_disk}) which may overlap themselves.
Consequently, some work is required to show that this actually defines a (not
necessarily regular) cell decomposition of the surface. Uniqueness, on the
other hand, is immediate from this definition.  We will also show that the
local Delaunay criterion implies the global empty circumcircle condition. A
Delaunay triangulation is obtained from \emph{the} Delaunay tessellation by
triangulating the non-triangular faces.  It follows that a Delaunay
triangulation, while in general not unique, differs from another Delaunay
triangulation only by edges with vanishing $\cot$-weights.  This is important
because it means that the discrete Laplace-Beltrami operator that will be
defined in Section~\ref{s.Laplace} depends only on the intrinsic geometry of
the surface.

\begin{definition}
\label{def:PF_surface}
A \emph{piecewise flat surface (PF surface)} $(S,d)$ is a 2-dimensional
differential manifold $S$, possibly with boundary, equipped with a metric $d$
which is flat except at isolated points, the \emph{cone points}, where $d$ has
cone-like singularities.
\end{definition} 

In other words, every interior point of a piecewise flat surface has a
neighborhood which is isometric to either a neighborhood of the Euclidean
plane or to a neighborhood of the apex of a Euclidean cone. The cone angle at
the apex may be greater than $2\pi$. (For a more detailed definition of
closed PF surfaces see Troyanov~\cite{Tro}.) In this paper, we consider only
\emph{compact surfaces} and we require the boundary (if there is a boundary)
to be piecewise geodesic. (The interior angle at a corner of the boundary may
be greater than $2\pi$.)

A \emph{tessellation}\/ of a PF surface is a cell decomposition such that the
faces are Euclidean polygons which are glued together along their edges. This
implies that the cone points and the corners of the boundary are vertices of
the cell decomposition. A \emph{triangulation}\/ is a tessellation where the
faces are triangles. For the following it is essential that we do
\emph{not}\/ require tessellations (and in particular triangulations) to be a
\emph{regular}\/ cell complexes, \ie~a glueing homomorphism may identify
points on the boundary of a cell.  For example, it is allowed that two edges
of a face may be glued to each other; and an edge may connect a vertex with
itself. \textit{A forteriori}, we do not require a tessellation to be
\emph{strongly regular}, \ie{} the intersection of two closed cells may not
be a single closed cell.

\begin{remark}
  From the intrinsic point of view, the carrier $S$ of a simplicial surface
  $\mathcal S$ with the metric induced by the ambient Euclidean space is a
  piecewise flat surface. The simplicial surface $\mathcal S$ also provides
  $S$ with a triangulation whose vertex set includes the cone points and the
  corners of the boundary. However, this triangulation is not intrinsically
  distinguished from other triangulations with the same vertex set.
\end{remark}

First, we will consider surfaces without boundary. To define the Delaunay
tessellation of a PF surface in terms of empty disks, we must allow an empty
disk to overlap with itself:

\begin{definition}[immersed empty disk]
  \label{def:immersed_empty_disk}
  Let $(S,d)$ be a compact PF surface \emph{without boundary}, and let
  $V\subset S$ be a finite set of points which contains all cone points. An
  \emph{immersed empty disk} is continuous map $\varphi:\bar D\rightarrow S$,
  where $D$ is an open round disk in the Euclidean plane and $\bar D$ is its
  closure, such that the restriction $\varphi|_{D}$ is an isometric immersion
  (\ie{} every $p\in D$ has a neighborhood which is mapped isometrically) and
  $\varphi(D)\cap V=\emptyset$.
\end{definition}

Hence any points in $\varphi^{-1}(V)$ are contained in the boundary of $D$:
$\varphi^{-1}(V)\subset\partial D$. 

\begin{definition}[Delaunay tessellation, no boundary]
  \label{def:delaunay_tesselation}
  Let $(S,d)$ be a compact PF surface \emph{without boundary}, and let
  $V\subset S$ be a finite set of points which contains all cone points. The
  \emph{Delaunay tessellation} of $(S,d)$ on the vertex set $V$ is the cell
  decomposition with the following cells: A subset $C\subset S$ is a closed
  cell of the Delaunay tessellation iff there exists an immersed empty disk
  $\varphi:\bar D\rightarrow S$ such that $\varphi^{-1}(V)$ is non-empty and
  $C$ is the image of the convex hull of $\varphi^{-1}(V)$:
  $C=\varphi(\conv\varphi^{-1}(V))$. The cell-attaching map is
  $\varphi|_{\conv\varphi^{-1}(V)}$; and the cell is a $0$-cell (vertex),
  $1$-cell (edge), or $2$-cell (face) if $\varphi^{-1}(V)$ contains one, two,
  or more points; respectively.
\end{definition}

\begin{proposition}
  The Delaunay tessellation as defined above is really a tessellation of
  $(S,d)$.
\end{proposition}

\begin{proof}
  Let us first remark that the vertex set of the Delaunay tessellation is
  obviously $V$. An edge $e$ is a geodesic segment such that there exists an
  immersed empty disk $\varphi:\bar D\rightarrow S$ with $\varphi^{-1}(V)$
  containing exactly two points, $\varphi^{-1}(V)=\{p_1,p_2\}$ such that
  $e=\varphi([p_1,p_2])$, where $[p_1,p_2]$ is the line segment joining $p_1$
  and $p_2$ in $D$. That the vertices and edges form a $1$-dimensional cell
  complex then follows from Lemma~\ref{lem:nocross} below. 
  
  Then we have to show that the open faces are indeed homeomorphic to open
  disks. (A cell attaching map $\varphi:|_{\conv\varphi^{-1}(V)}$ is
  \apriori{} only an immersion of the interior of the domain.) This follows
  from Lemma~\ref{lem:homeo}. 

  It is comparatively easy to see that the cell attaching homeomorphism
  $\varphi$ maps the boundary $\partial\conv\{p_1,p_2,\ldots,p_n\}$ into the
  1-skeleton; and that edges do not intersect open faces. We omit the details. 
  
  Finally, Lemma~\ref{lem:cellscoverS} asserts that every point in $S$ is
  contained in a closed cell.
\end{proof}

\begin{lemma}
  \label{lem:nocross}
  The edges do not cross each other or themselves.
\end{lemma}

\begin{proof}
  Let $e=\varphi([p_1,p_2])$ and $\tilde e=\tilde \varphi([\tilde p_1,\tilde
  p_2])$ be edges contained in the empty immersed disks $\varphi:\bar
  D\rightarrow S$ and $\tilde\varphi:\bar{\tilde D}\rightarrow S$,
  respectively, such that $p_1$, $p_2$ are the only points in
  $\varphi^{-1}(V)$ and $\tilde p_1$, $\tilde p_2$ are the only points in
  $\tilde\varphi^{-1}(V)$. Suppose $e$ and $\tilde e$ have an interior point
  in common: $\varphi(q)=\tilde\varphi(\tilde q)$ with $q\in(p_1,p_2)$ and
  $\tilde q\in(\tilde p_1,\tilde p_2)$. We are going to show that
  \begin{equation}
    \label{eq:e_equals_tilde_e}
    \varphi^{-1}(\tilde e)=[p_1,p_2].
  \end{equation}
  Since $\tilde\varphi(\tilde p_1),\tilde\varphi(\tilde
  p_2)\not\in\varphi(D)$ the Intersecting Chord Theorem implies
  \begin{equation*}
    \|p_1-q\|\cdot\|p_2-q\|
    \leq
    \|\tilde p_1-\tilde q\|\cdot\|\tilde p_1-\tilde q\|,
  \end{equation*}
  where $\|.\|$ denotes the Euclidean norm; see
  Figure~\ref{fig:chord-theorem}.
  \begin{figure}
    \centering
    \input{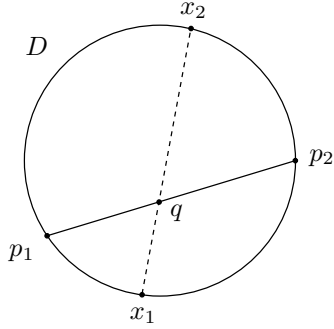}
    \caption{The Intersecting Chord
      Theorem says that $\|p_1-q\|\cdot\|p_2-q\|=\|x_1-q\|\cdot\|x_2-q\|$. }
    \label{fig:chord-theorem}   
  \end{figure}
  The opposite inequality follows equally. Hence
  $\varphi^{-1}(\tilde\varphi(\tilde p_1))$ and
  $\varphi^{-1}(\tilde\varphi(\tilde p_2))$ must be contained in $\partial
  D$. Since $p_1$ and $p_2$ are the only points in $\varphi^{-1}(V)$,
  this implies \eqref{eq:e_equals_tilde_e}.
\end{proof}

\begin{lemma}
  \label{lem:homeo}
  Let $\varphi: \bar D\rightarrow S$ be an immersed empty disk and suppose
  $\varphi^{-1}(V)=\{p_1,p_2,\ldots,p_n\}$ with $n\geq 3$. Let
  $P=\conv\{p_1,p_2,\ldots,p_n\}$. Then the restriction $\varphi|_{\interior
    P}$ of $\varphi$ to the interior of $P$ is injective (and hence a
  homeomorphism $\interior P\rightarrow\varphi(\interior P)$).
\end{lemma}

\begin{proof}
  Suppose $q\in \interior P$ and $\tilde q\in D$ with $q\not=\tilde q$\/ but
  $\varphi(q)=\varphi(\tilde q)$. We will show that $\tilde q\not\in\interior
  P$.
  
  Because $\varphi|_D$ is an isometric immersion, there is a neighborhood
  $q\ni U\subset D$ and an isometry $T$ of the Euclidean plane such that
  $T(q)=\tilde q$, $T(U)\subset D$, and $\varphi(T(x))=\varphi(x)$ for all
  $x\in U$.  Let $\tilde D=T(D)$ and $\tilde\varphi=\varphi\circ
  T^{-1}:\bar{\tilde D}\rightarrow S$. Since $\varphi$ and $\tilde\varphi$
  agree on the intersection $D\cap\tilde D$, there is a continuous map
  $\hat\varphi:\overline{D\cup\tilde D}\rightarrow S$ such that $\hat
  \varphi|_{\bar D}=\varphi$ and $\hat\varphi|_{\bar{\tilde D}}=\tilde
  \varphi$.
  
  \begin{figure}
    \centering
    \input{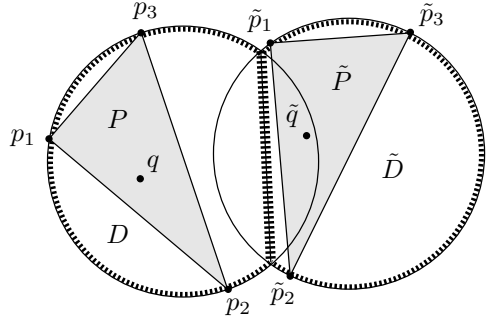}
    \caption{The dashed lines mark the boundaries of of the circular segments 
      $\conv((\partial D)\setminus\tilde D)$ and $\conv((\partial\tilde
      D)\setminus D)$.}
    \label{fig:facesaredisks}
  \end{figure}
  
  Now let $\tilde p_i=T(p_i)$ ($i=1,\ldots n$) and $\tilde
  P=T(P)=\conv\{\tilde p_1,\tilde p_2,\ldots,\tilde p_n\}$. Then $P$ and
  $\tilde P$ have no common interior points: $\interior P\cap\interior\tilde
  P=\emptyset$. Indeed, since $D$ and $\tilde D$ are ``empty'',
  $\{p_1,p_2,\ldots,p_n\}\subset(\partial D)\setminus\tilde D$ and $\{\tilde
  p_1,\tilde p_2,\ldots,\tilde p_n\}\subset(\partial\tilde D)\setminus D$.
  Hence $P\subset\conv((\partial D)\setminus\tilde D)$ and $\tilde
  P\subset\conv((\partial\tilde D)\setminus D)$. But the circular segments
  $\conv((\partial D)\setminus\tilde D)$ and $\conv((\partial\tilde
  D)\setminus D)$ have no interior points in common. (See
  Figure~\ref{fig:facesaredisks}.)
  
  Now $\tilde q\in\interior\tilde P$ implies $\tilde q\not\in\interior P$.
\end{proof}

\begin{lemma}
  \label{lem:cellscoverS}
  Every point $x\in S$ is contained in a closed cell.
\end{lemma}

\begin{proof}
  Consider two immersed disks $\varphi:\bar D\rightarrow S$,
  $\tilde\varphi:\tilde{\bar D}\rightarrow S$ as equivalent if they differ
  only by a change of parameter, \ie\ if $\tilde\varphi=\varphi\circ T$ for
  some isometry of the Euclidean plane. The manifold of equivalence classes
  is parameterized by the set
  \begin{equation*}
    \mathcal D =\big\{(c,r)\in S\times\R_{>0}\,\big|\,d(c,V)\geq r\big\}
  \end{equation*}
  of center/radius pairs.  If we adjoin degenerated immersed empty disks with
  radius $0$, we obtain a compact manifold with boundary, parameterized by
  \begin{equation*}
    \bar{\mathcal D} 
    =\big\{(c,r)\in S\times\R_{\geq 0}\,\big|\,d(c,V)\geq r\big\}.
  \end{equation*}
  For a point $x\in S$ the \emph{power function} with respect to $x$ is the
  continuous function 
  \begin{equation*}
    \power_x:\bar{\mathcal D}\rightarrow\R,\quad
    \power_x(c,r)=\big(d(c,x)\big)^2-r^2.
  \end{equation*}
  If $\varphi:\bar D\rightarrow S$ is an immersed empty disk with center $c$
  and radius $r$, then $\power_x(c,r)$ is smaller than, equal to, or greater
  than zero depending on whether $x\in\varphi(D)$, $x\in\partial\varphi(D)$,
  or $x\in S\setminus\varphi(\bar D)$.
  
  Let $x\in S$. We have to show that $x$ is contained in some closed cell. If
  $x\in V$ this is clear, because the points in $V$ are $0$-cells. So assume
  $x\in S\setminus V$. Since $\power_x$ is continuous on the compact set
  $\bar{\mathcal D}$, there is a
  $(c_{\text{min}},r_{\text{min}})\in\bar{\mathcal D}$ where $\power_x$
  attains its minimum. Since $x\not\in V$, there is an empty disk containing
  $x$ and hence $\power_x(c_{\text{min}},r_{\text{min}})<0$. Let
  $\varphi:\bar D\rightarrow S$ be an immersed empty disk with center
  $c_{\text{min}}$ and radius $r_{\text{min}}$; \ie~$D$ is a disk in $\R^2$
  with center $m\in\R^2$ and radius $r_{\text{min}}$ and
  $\varphi(m)=c_{\text{min}}$. There is a $p\in D$ with $\varphi(p)=x$ and
  $\|p-m\|=d(x,c_{\text{min}})$. We show by contradiction that
  $p\in\conv\varphi^{-1}(V)$.
  
  Suppose the opposite is true: $p\not\in\conv\varphi^{-1}(V)$. Then there
  exists a closed half-space $H\subset\R^2$ with $\varphi^{-1}(V)\subset H$
  but $p\not\in H$. In that case there exists another immersed empty disk
  $\tilde\varphi:\tilde D\rightarrow M$ with $D\setminus\tilde
  D\subset\interior H$, $\tilde D\setminus D\subset\R^2\setminus H$, and
  $\varphi|_{D\cap\tilde D}=\tilde\varphi|_{D\cap\tilde D}$ (see
  Figure~\ref{fig:smallerpower}). 
  \begin{figure}
    \centering
    \input{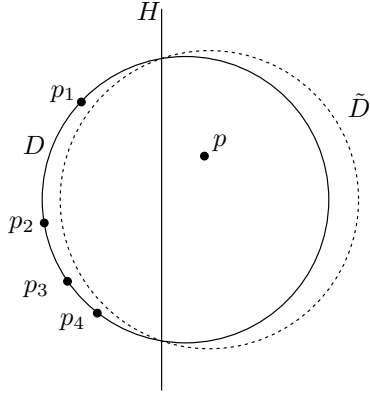}
    \caption{If $p$ is not contained in the convex hull of
      $\varphi^{-1}(V)=\{p_1,\ldots,p_n\}$, then there is another immersed
      empty disk with smaller power. (See proof of
      Lemma~\ref{lem:cellscoverS}.)}
    \label{fig:smallerpower}
  \end{figure}
  Let $\tilde m$ and $\tilde r$ be center and radius of
  $\tilde D$. Then
  \begin{equation*}
    \|p-\tilde m\|^2-{\tilde r}^2 < \|p-m\|^2-r^2.
  \end{equation*}
  (To see this, note that $q\mapsto (\|q-\tilde m\|^2-{\tilde
    r}^2)-(\|q-m\|^2-r^2)$ is an affine linear function $\R^2\rightarrow\R$,
  which vanishes on $\partial H$ and is positive on $\interior H$, and
  negative on $\R^2\setminus H$.)  This implies $\power_x(\tilde c,\tilde
  r)<\power_x(c_{\text{min}},r_{\text{min}})$, where $\tilde
  c=\tilde\varphi(\tilde m)$. This contradicts the assumption that $\power_x$
  attains its minimum on $(c_{\text{min}},r_{\text{min}})$.
\end{proof}

\paragraph{Delaunay tessellations of PF surfaces with boundary.} 

Now let $(S,d)$ be a compact PF manifold with piecewise geodesic boundary.
Let $V\subset S$ be a finite set of points which contains all cone points and
all corners of the boundary. The boundary is then the union of geodesic
segments connecting points in $V$ but containing no points of $V$ in their
interior.  Each connected component of the boundary is a closed geodesic
polygon with vertices in $V$. To each such boundary polygon glue a PF surface
obtained by cyclicly gluing together the appropriate number of isosceles
triangles with appropriate base lengths and legs of length $R>0$. Each of
these caps contains a special point where the triangle apices are glued
together. It is in general a cone point. The result of closing all wholes
with such caps is a closed PF surface $(\hat S, d)$. Let $\hat V$ be the
union of $V$ and the set of special points of the caps. If $R$ is chosen
large enough, then the isosceles triangles in the caps will be faces of the
Delaunay tessellation of $(\hat S, d)$ with respect to $\hat V$. (This is so
because if $R$ is large enough, the immersed circumdisks intersect $S$ in
lunes which are so small that they are empty.) The Delaunay tessellation of
the bounded surface $(S,d)$ with respect to $V$ is defined to be the cell
complex obtained by removing these triangles.

\paragraph{Delaunay triangulations.} A \emph{Delaunay triangulation} is a
triangulation obtained from a Delaunay tessellation by triangulating all
non-triangular faces. A Delaunay triangulation is characterized by the
\emph{empty circumcircle property:} A tessellation of $(S,d)$ on the vertex
set $V$ is a Delaunay triangulation iff for each face $f$ there exists an
immersed empty disk $\varphi:\bar D\rightarrow S$ such that there are three
points $p_1,p_2,p_3\in\varphi^{-1}(V)$ with
$f=\varphi(\conv\{p_1,p_2,p_3\}$). (But there may be more than three points in
$\varphi^{-1}(V)$.)

In Proposition~\ref{t.local->global_Delaunay} we give a more local
characterization of Delaunay triangulations. 

\begin{definition}
  \label{def:locally_delaunay}
  Let $T$ be a geodesic triangulation of $(S,d)$ with vertex set $V$, and let
  $e$ be an interior edge of $T$. Since all faces of $T$ are isometric to
  Euclidean triangles we can isometrically unfold to the plane the two
  triangles of $T$ that are adjacent to $e$. We say that $e$ is \emph{locally
    Delaunay} if the vertices of the two unfolded triangles are not contained
  inside the circumcircles of these triangles.
\end{definition}

For our investigation of discrete Laplace operators we will need the
following characterization of Delaunay edges.

\begin{lemma}
  \label{l.Delaunay_edge}
  An interior edge $e$ of a triangulation $T$ of a piecewise flat surface
  $(S,d)$ is locally Delaunay if and only if the sum of the angles opposite
  $e$ in the adjacent triangles does not exceed $\pi$.
\end{lemma}

This follows immediately from the fact that opposite angles in a
circular quadrilateral sum to $\pi$.

Clearly all interior edges of a Delaunay triangulation are locally Delaunay.
This property actually characterizes Delaunay triangulations:

\begin{proposition}
  \label{t.local->global_Delaunay}
  A triangulation $T$ of a piecewise flat surface $(S,d)$ is a Delaunay
  triangulation if and only if all interior edges of $T$ are locally Delaunay.
\end{proposition}

The following proof is an adaptation of Delaunay's original
argument~\cite{Del} for Delaunay triangulations in $\R^n$. (See also
Edelsbrunner~\cite[p.~8]{edel} for a more easily available modern
exposition.)  Alternatively, one could also adapt the argument of Aurenhammer
and Klein~\cite{aurenhammer00} for Delaunay triangulations in the plane.

\begin{proof}
  If $S$ is a manifold with boundary, construct a closed PF surface by
  glueing piecewise flat disks to the boundary components in the manner
  described above. If $R$ is chosen large enough, all edges will be locally
  Delaunay. It remains to prove the Proposition for closed PF surfaces.
  
  Suppose that all interior edges of the triangulation $T$ of the closed PF
  surface $(S,d)$ are locally Delaunay. We want to show that $T$ is a
  Delaunay triangulation. To this end we will show that the empty
  circumcircle property holds. Let $t_0$ be a triangle of $T$.  Starting with
  $t_0$ we develop a part of $T$ in the Euclidean plane (see
  Figure~\ref{fig:layout}, left).
  \begin{figure}
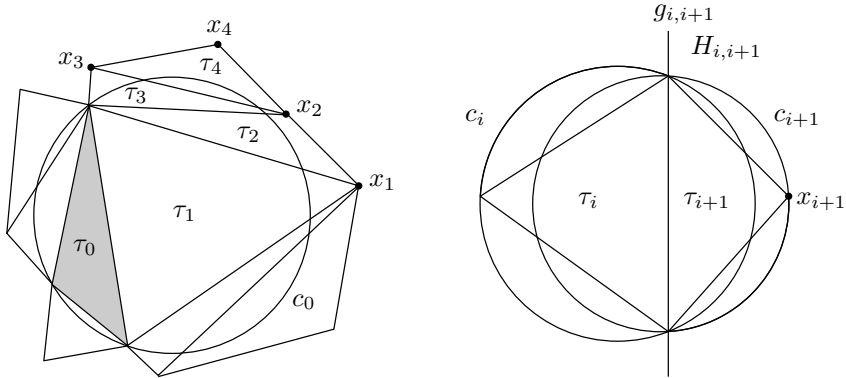

    \hfill \input{layout.tex} \hfill \input{tausequence.tex}
    \hspace*{\fill}
    \caption{\emph{Left:}~The layed out triangles. \emph{Right:}~The power line
      of $c_i$ and $c_{i+1}$ is $g_{i,i+1}$. Since the edge $e_{i,i+1}$ is
      assumed to be Delaunay, $p_{c_i}(x)\geq p_{c_{i+1}}(x)$ in
      $H_{i,i+1}$.}
    \label{fig:layout}
  \end{figure}
  Begin with a triangle $\tau_0$ in the Euclidean plane that is congruent to
  $t_0$. Let $c_0$ be the circumcircle of $\tau_0$.  Next, lay out congruent
  copies of the triangles neighboring $t_0$.  This introduces new vertices in
  the plane, which are on or outside $c_0$ because the edges of $t_0$ are
  Delaunay by assumption.  Keep laying out neighboring triangles in the plane
  at free edges but only if the free edges intersect $c_0$. (Different layed
  out triangles may correspond to the same triangle in $(S,d)$.) Each new
  triangle introduces a new vertex and we will show that they do not lie
  inside $c_0$.  Hence, when the layout process stops (when all free edges do
  not intersect $c_0$), the triangles simply cover the inside of $c_0$. It
  follows that $t_0$ has an empty circumcircle in $(S,d)$.
  
  Let $\tau_0, \tau_1, \ldots,\tau_n$, be a sequence of layed out triangles
  such that $\tau_i$ and $\tau_{i+1}$ share an edge $e_{i,i+1}$. Let
  $x_{i+1}$ be the vertex opposite $e_{i,i+1}$ in $\tau_{i+1}$. Assuming that
  $x_i$ is not inside $c_0$ for all $i<n$ we will show that the same holds
  for $x_n$.  Let $g_{i,i+1}$ be the straight line containing $e_{i,i+1}$ and
  let $H_{i,i+1}$ be the half space bounded by $g_{i,i+1}$ and containing
  $\tau_{i+1}$.  Then
  \begin{equation*}
    (H_{0,1}\cap D_0)\supset(H_{1,2}\cap D_0)\supset\ldots
    \supset(H_{n-1,n}\cap D_0),
  \end{equation*}
  where $D_0$ is the open disk bounded by $c_0$. Hence it remains to consider
  the case where $x_n\in H_{i,i+1}$ for all $i=0,\ldots,n-1$, because
  otherwise $x_n\not\in D_0$.
  
  Now consider the power of a point $x\in\R^2$ with respect to a circle
  $c\subset\R^2$ with center $x_c$ and radius $r$ as a function of $x$:
  \begin{equation*}
    p_c(x)=\|x-x_c\|^2 - r^2.
  \end{equation*}
  It is positive, zero, or negative if $x$ lies outside, on, or inside $c$,
  respectively. The \emph{power line}\/ of two different circles $c$ and $c'$
  is the locus of points $x$ with $p_c(x)=p_{c'}(x)$. It is a straight line
  because $p_c(x)-p_{c'}(x)$ is linear in $x$. The power line of two
  intersecting circles is the line through the intersection points. Let $c_i$
  be the circumcircle of $\tau_i$. Either $c_i=c_{i+1}$ or the power line of
  $c_i$ and $c_{i+1}$ is $g_{i,i+1}$ and
  \begin{equation}
    \label{e.H_and_p}
    H_{i,i+1}=\{x:p_{c_i}(x)\geq p_{c_{i+1}}(x)\}.
  \end{equation}
  (See Figure~\ref{fig:layout}, right.) Indeed, $p_{c_{i+1}}(x_{i+1})=0$ and
  since the edge $e_{i,i+1}$ is locally Delaunay by assumption,
  $p_{c_i}(x_{i+1})\geq 0$. Hence
  \begin{equation*}
    p_{c_0}(x_n)\geq p_{c_1}(x_n)\geq\ldots\geq p_{c_n}(x_n)=0,
  \end{equation*}
  and therefore $x_n$ lies on or outside $c_0$. This concludes the proof.
\end{proof}

The \emph{edge flipping algorithm}\/ may be used to construct a Delaunay
triangulation of a piecewise flat surface $(S,d)$ with marked points
$V\subset S$:
\begin{enumerate}\setlength{\itemsep}{0pt}
\item Start with any triangulation $T$ of $(S,d)$ with vertex set
  $V$.
\item\label{item_repeat} If all interior edges of $T$ are locally Delaunay,
  stop.
\item Otherwise there is an interior edge $e$ of $T$ which is not locally
  Delaunay. Perform an \emph{intrinsic edge flip}: Replace $e$ by the
  other diagonal of the quadrilateral formed by the two triangles adjacent to
  $e$.  Go to Step~\ref{item_repeat}.
\end{enumerate}
The following two propositions show that this is indeed an algorithm.
\begin{proposition}
  \label{prop:flippable}
  If an edge is not locally Delaunay, then it can be flipped.
\end{proposition}
\begin{proof}
  In $\R^2$, an edge is flipable iff the two adjacent triangles form a
  convex quadrilateral. In a PF surface, there is an additional topological
  obstruction to flipability: Since the triangulation may not be regular, an
  edge may be adjacent to the same triangle on both sides. So suppose the
  edge $e$ is not locally Delaunay, \ie{} the sum of opposite angles exceeds
  $\pi$. Then there are two different triangles adjacent to $e$, because the
  sum of all angles in one triangle is $\pi$. These two triangles form a
  Euclidean quadrilateral (possibly with some of the boundary edges
  identified with each other), which is convex by the usual argument. Hence
  $e$ can be flipped.
\end{proof}
\begin{proposition}[Indermitte et~al.~\cite{ILTC}]
  \label{t.flipalgorithm}
  The edge flipping algorithm terminates after a finite number of steps.
\end{proposition}
Together with Proposition~\ref{t.local->global_Delaunay} this implies that
the edge flipping algorithm produces a Delaunay triangulation (in the global
empty-circumcircle sense). To prove Proposition~\ref{t.flipalgorithm}, one
has to show that the algorithm cannot loop infinitely. In the setting of
\emph{planar} Delaunay triangulations, it is enough to define a suitable real
valued function on the set of triangulations on the given vertices which
decreases (or increases) with each edge flip. Because this set of
triangulations is finite, the algorithm has to terminate. As a further
consequence such a function attains its minimal (or maximal) value on the
Delaunay triangulations. Several such functions are known, see for example
Lambert~\cite{lambert}, Musin~\cite{Mu}, Rivin~\cite[Sec.~10]{Riv}, and the
survey article~\cite{aurenhammer00}. When we consider Delaunay triangulations
of PF surfaces, however, the set of triangulations on the marked points may
be infinite. (The fact that the number of combinatorial types of
triangulations is finite~\cite{Riv} is not sufficient to make the argument.)
For example, the surface of a cube has infinitely many geodesic
triangulations on the eight vertices. To prove
Proposition~\ref{t.flipalgorithm} by means of a function which decreases with
every flip, one has to show that it has the following additional property.

\begin{definition}
  Let $\mathcal T$ be the set of triangulations of a PF surface on a given
  set of marked points and let $f:\mathcal T\rightarrow\R$. We say that $f$
  is \emph{proper} if for any $M\in\R$ the number of triangulations
  $T\in\mathcal T$ with $f(T)\leq M$ is finite.
\end{definition}

In their proof of Proposition~\ref{t.flipalgorithm}, Indermitte
\textit{et~al.}~\cite{ILTC} use the sum of squared circumcircle radii as
proper function which decreases with every flip. (For a proof that it
decreases with every flip, they refer to an unpublished PhD thesis. However,
see Musin~\cite{Mu}, his Theorem~3 and the following Lemma, for a hint on how
to prove this.) Another possible choice is Musin's harmonic index. Below we
show that it is proper and decreases with every flip. The latter fact we
deduce from Rippa's Theorem~\cite{R}, which is also of independent interest
in connection with the Laplace-Beltrami operator.  Rippa's proof holds
without change also for piecewise flat surfaces.

\begin{rippastheorem}%[Rippa~\cite{R}]
%\label{t.Dirichlet-Delaunay}
Let $(S,d)$ be a piecewise flat surface and let $V\subset S$ be a set of
marked points which contains the cone points and the corners of the boundary.
Let $f:V\to {\R}$ be a function on the marked points. For each triangulation
$T$ of $(S,d)$ with vertex set $V$ let $f_T:S\to\R$ be the PL interpolation
of $f$ that is linear on the faces of the triangulation~$T$.

Suppose $T_1$ is a triangulation with an interior edge $e$ and $T_2$ is
obtained from $T_1$ by flipping $e$. If the edge $e$ is a Delaunay edge after
the flip, \ie~in $T_2$, then
\begin{equation*}
  E(f_{T_1})\geq E(f_{T_2}),
\end{equation*}
where $E$ denotes the Dirichlet energy as in Section~\ref{s.energy}.
Equality holds only if $f_{T_1}=f_{T_2}$ or if $e$ was also a Delaunay edge
in $T_1$.

As a consequence, the minimum of the Dirichlet energy among all possible
triangulations is attained on the Delaunay triangulations $(S,d)$:
$$
\min_T \int_S \mid \nabla f_T\mid^2=\int_S \mid \nabla
f_{T_D}\mid^2,
$$
where $T_D$ is any Delaunay triangulation.

Moreover, for generic $f:V\to {\R}$, this property of Delaunay
triangulations is characteristic.
\end{rippastheorem}

Rippa's proof is based on the following comparison formula for the
Dirichlet energies of two possible triangulations of a
quadrilateral
$$
 E(f_{T_1})-E(f_{T_2})=\frac{(f_1-f_2)^2}{2\sin \theta} \frac{(r_1+
 r_3)(r_2+r_4)}{r_1 r_2 r_3 r_4}(r_1 r_3- r_2 r_4).
$$
Here $T_1$ and $T_2$ are the two triangulations of the convex
quadrilateral $Q=(x_1,x_2,x_3,x_4)$ obtained by addition of the
diagonals $(x_1,x_3)$ and $(x_2,x_4)$ respectively, $f_1$ and
$f_2$  are the values of $f_{T_1}$ and $f_{T_2}$ at the
intersection point $x_0$ of the diagonals, $r_1,\ldots,r_4$ are
the distances from $x_0$ to the vertices $x_1,\ldots,x_4$ of the
quadrilateral and $\theta$ is the intersection angle of the
diagonals.

For a Euclidean triangle $t$ with sides $a$, $b$, $c$ and area $A$,
Musin~\cite{Mu} defines the \emph{harmonic index}\/ as
\begin{equation*}
\hrm(t)=\frac{a^2+b^2+c^2}{A}.
\end{equation*}
The harmonic index of a triangulation $T$ with face set $F$ is the sum of the
harmonic indices of all triangles:
\begin{equation*}
  \hrm(T)=\sum_{t\in F}\hrm(t).
\end{equation*}

\begin{proposition}
  The harmonic index is proper.
\end{proposition}

\begin{proof}
  In a PF surface, there may be infinitely many geodesic lines connecting two
  points, but for any $L\in\R$ only a finite number of them have length $\leq
  L$~\cite{ILTC}. Hence there are only a finite number of triangulations of a
  PF surface with marked points such that all edges are not longer than $L$.
  Now the Proposition follows from the inequality
  \begin{equation*}
    hrm(T)\geq \frac{l_{\text{max}}(T)}{A_{\text{tot}}},
  \end{equation*}
  where $l_{\text{max}}(T)$ is the largest length of an edge of $T$ and
  $A_{\text{tot}}$ is the total area of the PF surface. Indeed, if
  $hrm(T)\leq M$, then $l_{\text{max}}(T)\leq M A_{\text{tot}}$, and there
  are only finitely many triangulations satisfying this bound on edge lengths.
\end{proof}

The following theorem was stated by Musin without proof~\cite{Mu}.
\begin{theorem}
\label{t.musin}
  With the notation and under the conditions of
  Rippa's Theorem%~\ref{t.Dirichlet-Delaunay},
  \begin{equation*}
    \hrm(T_1)\geq\hrm(T_2)
  \end{equation*}
  and equality holds only if $e$ is a Delaunay edge in $T_1$ as well. This
  implies that the harmonic index is minimal for a Delaunay triangulation
  $T_D$ (and hence for all of them):
  \begin{equation*}
    \min_T \hrm(T) = \hrm(T_D).
  \end{equation*}
\end{theorem}
\begin{proof}
  The harmonic index of a triangle is
  \begin{equation*}
    \hrm(t)=4(\cot\alpha+\cot\beta+\cot\gamma),
  \end{equation*}
  where $\alpha$, $\beta$, $\gamma$ are the angles of the triangle.
  (Because the area is $A=\frac{1}{2}a h_a=\frac{1}{2}b h_b=\frac{1}{2}c
  h_c$, where $h_a, h_b, h_c$ are the heights of the triangle; and
  $a/h_a=\cot\beta+\cot\gamma$, \emph{etc.})

  For a triangulation $T$ of $(S,d)$ with vertex set $V$ and $x\in
  V$ let $\delta_{x,T}:S\to\R$ be the function that is linear on the
  triangles of $T$, equal to $1$ at $x$, and equal to $0$ at all other marked
  points in $V$. Then
  \begin{equation*}
    \sum_{x\in V}E(\delta_{x,T})=\frac{1}{2}\sum_{\text{angles $\alpha$ in $T$}}
    \cot\alpha=\frac{1}{8}\hrm(T).
  \end{equation*}
  Hence the theorem follows from Rippa's Theorem.%~\ref{t.Dirichlet-Delaunay}.
\end{proof}

\section{The discrete Laplace-Beltrami operator and discrete harmonic
  functions}
\label{s.Laplace}

We are now in a position to define the discrete Laplace operator on a
simplicial surface in an intrinsic way.
\begin{definition}
  \label{d.Laplace}
  Let $\mathcal S$ be a simplicial surface with vertex set $V$ and let $S$ be
  its carrier, which is a piecewise flat surface.  The \emph{discrete
    Laplace-Beltrami operator $\Delta$ of a simplicial surface} $\mathcal S$
  is defined as follows. For a function $f:V\to\R^n$ on the vertices, the
  value of $\Delta f: V\to\R^n$ at $x_i\in V$ is
  \begin{equation}                            \label{eq.Laplace_D}
    \Delta f(x_i)=\sum_{x_j\in V:(x_i,x_j)\in E_D}
    \nu(x_i,x_j)(f(x_i)-f(x_j)),
  \end{equation}
  where $E_D$ is the edge set of a Delaunay triangulation of $S$ and the
  weights are given by
  \begin{equation}
  \label{eq:weights}
  \nu(x_i,x_j)=
  \left\{
    \begin{array}{ll}
      \frac{1}{2}(\cot\alpha_{ij}+\cot\alpha_{ji})
      & \text{for interior edges}\\
      \frac{1}{2}\cot\alpha_{ij} & \text{for boundary edges}
    \end{array}
  \right..
 \end{equation}
 Here $\alpha_{ij}$ (and $\alpha_{ji}$ for interior edges) are the angles
 opposite the edge $(x_i,x_j)$ in the adjacent triangles of the Delaunay
 triangulation (see Figure~\ref{f.angles}).

  The \emph{discrete Dirichlet energy of $f$} is
  $$
  \mathcal{E}_D=\frac{1}{2}\sum_{(x_i,x_j)\in E_D}
  \nu(x_i,x_j)(f(x_i)-f(x_j))^2.
  $$
\end{definition}

Due to Lemma \ref{l.Delaunay_edge} and the formula $
\cot\alpha+\cot\beta=\frac{\sin(\alpha + \beta)}{\sin\alpha \sin\beta}$ the
discrete Laplace operator has non-negative weights. The edges with vanishing
weights are diagonals of non-triangular cells of the Delaunay tessellation.
Erasing such edges in (\ref{eq.Laplace_D}) we obtain a discrete Laplace
operator on the Delaunay tessellation of $S$ with positive weights
$\nu(x,x_i)$. Moreover, this property is characteristic for Delaunay
triangulations:  Consider a piecewise flat surface $(S,d)$ with a
triangulation $T$. Denote by $\Delta_T$ the Laplace operator of the
triangulation $T$: it is given by the same formula (\ref{eq.Laplace_D}) with
the weights $\nu_T$ determined by the triangulation $T$ by the same formulas
(\ref{eq:weights}) as for the Delaunay triangulation.
The following observation is elementary.
\begin{proposition}
  The Laplace operator $\Delta_T$ of the triangulation $T$ has non-negative
  weights $\nu_T$ if and only if the triangulation $T$ is Delaunay.
\end{proposition}

Laplace operators with positive weights on graphs possess
properties analogous to the smooth theory.

\begin{definition}
  A discrete function $f:V\to{\R^n}$ on a simplicial surface is
  \emph{harmonic}\/ if $\Delta f(x)=0$ for all interior vertices $x$.
\end{definition}
Discrete harmonic functions satisfy the maximum principle: A real valued
harmonic function attains its maximum on the boundary. This implies:
\begin{proposition}
  \label{p.convex_hull}
  For each interior vertex $x$, the value $f(x)$ of a harmonic function
  $f:V\to {\R}^n$ lies in the convex hull of the values $f(x_i)$ on its
  neighbors.
\end{proposition}

We conclude this section with some standard facts regarding boundary value
problems for the discrete Laplace operator. Let us fix a subset
$V_{\partial}\subset V$ of vertices--which may but need not be the set of
boundary vertices--and a function $g:V_{\partial}\to{\R}$. The problem of
finding a function that satisfies $\Delta f (x)=0$ for $x\in V\setminus
V_\partial$ and $f(x)=g(x)$ for $x\in V_\partial$ is called a \emph{Dirichlet
  boundary value problem}.  The problem of finding a function that satisfies
$\Delta f (x)=0$ for $x\in V\setminus V_\partial$ and $\Delta f(x)=g(x)$ for
$x\in V_\partial$ is a \emph{Neumann boundary value problem}.

\begin{theorem}
  For arbitrary $V_\partial$ and $g(x)$ the solutions of the Dirichlet and
  Neumann boundary value problems exist and are unique (up to an additive
  constant if $V_\partial=\emptyset$ and in the Neumann case). The solution
  minimizes the Dirichlet energy $\mathcal{E}_D$.
\end{theorem}

Indeed, solutions of the Dirichlet and Neumann boundary value problems are
critical points of the Dirichlet energy. Since the energy is a positive
definite quadratic form, its only critical point is the global minimum.

More complicated boundary conditions, such as so called ``natural
boundary conditions'' (see Desbrun et~al.~\cite{DMA}), are also
intensively used in geometry processing. The corresponding
existence and uniqueness results are still missing.

\section{Discrete mean curvature and minimal surfaces}
\label{s.minimal}

In this section we adapt the definitions of the mean curvature
vector for simplicial surfaces and minimal surfaces originally
suggested by Pinkall and Polthier \cite{PP} to the discrete
Laplace operator that we introduced in Section~\ref{s.Laplace}.

For a \emph{smooth}\/ immersed surface $f:\R^2\supset U\to{\R}^3$ the mean
curvature vector is given by the formula $H=\Delta f$, where $\Delta$ is the
Laplace-Beltrami operator of the surface. For a simplicial surface we define
the mean curvature vector by the same formula, following~\cite{PP}, but we
use a different Laplace operator:
\begin{definition}
  \label{d.H1}
  Let $\mathcal S$ be a simplicial surface with carrier $S$. The
  \emph{discrete mean curvature vector} at a vertex $x$ is
  \begin{equation*}
    {\mathcal H}(x)= \Delta f(x),
  \end{equation*}
  where $f:S\to\R^3$ is the restriction of the identity map on $\R^3$ to $S$,
  and $\Delta$ is the discrete Laplace-Beltrami operator of
  Definition~\ref{d.Laplace}.
\end{definition}

The discrete mean curvature vector $\mathcal H(x)$ of a simplicial
surface corresponds to the integral of the mean curvature vector
of a smooth surface over a neighborhood of the point $x$. (Note
when we scale the simplicial surface, $\mathcal H$ varies as the
integral of the mean curvature over some domain.) This suggests
the following alternative definition:
\begin{definition}
  \label{d.H2}
  In the setup of Definition~\ref{d.H1}, let $C(x)$ be a Voronoi cell of the
  vertex $x$ of a simplicial surface. The \emph{discrete mean curvature
    vector density} at~$x$  is defined by
  $$
  H(x)=\frac{\Delta f(x)}{A(C(x))},
  $$
  where $\Delta$ is the discrete Laplace-Beltrami operator of
  Definition~\ref{d.Laplace} and $A(C(x))$ is the area of the Voronoi cell
  $C(x)$.
\end{definition}

This definition is similar to the definition of mean curvature
suggested by Meyer, Desbrun, Schr\"oder and Barr~\cite{MDSB}.
Again, the difference is that we contend that one should use the
discrete Laplace-Beltrami operator.

\begin{definition}[Wide definition of simplicial minimal surfaces]
\label{d.minimal} A simplicial surface is called \emph{minimal}\/
if its mean curvature vector vanishes identically.
\end{definition}

Since the embedding $f:S\to \R^3$ of a simplicial minimal surface is
harmonic, Proposition~\ref{p.convex_hull} implies

\begin{proposition}
Every interior vertex of a simplicial minimal surface lies in the convex
hull of its neighbors.
\end{proposition}

The following stricter definition is also natural.

\begin{definition}[Narrow definition of simplicial minimal surfaces]
  A simplicial surface $\mathcal S$ is called \emph{minimal}\/ if its mean
  curvature vector vanishes identically and the intrinsic Delaunay
  triangulation of the carrier of $\mathcal S$ coincides with the
  triangulation induced by the simplicial complex $\mathcal S$.
\end{definition}

%% It can happen that the (intrinsic) Delaunay triangulation of the carrier~$S$
%% of a simplicial minimal surface~$\mathcal S$ (by Definition~\ref{d.minimal})
%% coincides with the triangulation induced by the simplicial complex $\mathcal
%% S$. In that case, $\mathcal S$ is also a discrete minimal surface in the
%% sense of Pinkall and Polthier \cite{PP}. Then, the surface 
Such a simplicial minimal surface is a critical point of the area functional
under variations of the vertex positions~\cite{PP}. We would like to note
that there exists also a non-linear theory of discrete minimal surfaces based
on the theory of circle patterns~\cite{BHS}.

The mean curvature flow for simplicial surfaces is given by the
equation
$$
\frac{df}{dt}(x)=H(x).
$$
This flow may change the Delaunay triangulation of the surface. At some
moment two Delaunay circles coincide and the diagonal of the corresponding
quadrilateral flips. However, since the weights of the diagonals vanish at
the flip moment, the discrete Laplace-Beltrami operator and mean curvature
vectors are continuous functions of time $t$.

For the numerical computation of discrete minimal surfaces one should use the
algorithm of~\cite{PP} with an extra step to adapt the Delaunay triangulation:
Start with a simplicial surface $\mathcal S_0$ which respects the
given boundary conditions. Calculate the Delaunay triangulation of
the carrier $S_0$ and the weights $\nu$. Find an $f:S_0\to\R^3$
which respects the boundary conditions and minimizes the Dirichlet
energy (see Definition~\ref{d.Laplace}). You may start with the
embedding of $S_0$ as initial guess. You obtain a new simplicial
surface $\mathcal S_1$ which is combinatorially equivalent to
$\mathcal S_0$ but geometrically different. Calculate the Delaunay
triangulation and weights $\nu$ of $S_1$ and find an
$f:S_1\to\R^3$ that minimizes the Dirichlet energy. Iterate.

\paragraph{Acknowledgments.} We would like to thank Peter Schr\"oder for
useful discussions and for pointing out Rippa's paper to us. G\"unter~M.
Ziegler and Ivan Ismestiev have made valuable comments for which we are
grateful. Last but not least we would like to thank the referees for their
helpful suggestions, in particular for informing us of the
reference~\cite{duffin1}.

%--------------------------------------------------------------------------

\end{document}